\documentclass[12pt]{article}
\usepackage{amssymb
}
\usepackage{amsmath}
\allowdisplaybreaks

\newcommand{\rom}[1]{{\rm #1}}

\newtheorem{theorem}{Theorem}

\newtheorem{lemma}{Lemma}
\newtheorem{remark}{Remark}
\newtheorem{corollary}{Corollary}

\begin{document}

\begin{center}{\Large \bf
 On the correlation measure of a family of commuting Hermitian operators
with applications to particle densities of the quasi-free representations of the CAR and CCR}\end{center}


{\large Eugene Lytvynov}\\
Department of Mathematics, University of Wales Swansea, Singleton
Park, Swansea SA2 8PP, U.K.\\ e-mail:
\texttt{e.lytvynov@swansea.ac.uk}\vspace{2mm}

{\large Lin Mei}\\
Department of Mathematics, University of Wales Swansea, Singleton
Park, Swansea SA2 8PP, U.K.\\ e-mail:
\texttt{300957@swansea.ac.uk}\vspace{2mm}

{\small
\begin{center}
{\bf Abstract}
\end{center}
\noindent
Let $X$ be a   locally compact, second countable Hausdorff
topological space.
We consider a family of commuting Hermitian operators $a(\Delta)$ indexed
by all measurable, relatively compact sets $\Delta$ in $X$ (a quantum stochastic process over $X$). For such a family, we introduce the notion of a correlation measure. We prove that, if the family of operators possesses a correlation measure which satisfies some condition of growth, then there exists
a point process over $X$ having the same correlation measure. Furthermore, the operators
$a(\Delta)$ can be realized as multiplication operators in the $L^2$-space with respect to this point process.
In the proof, we utilize the notion of $\star$-positive definiteness, proposed in
\cite{Kuna1}. In particular, our result extends the criterion of existence of a point process from that paper to the case
of the topological  space $X$, which is a standard underlying space in the theory of point processes.
As applications, we discuss particle densities of the quasi-free representations of the CAR and CCR, which lead
to fermion, boson,  fermion-like, and boson-like (e.g.\ para-fermions and para-bosons of order 2) point processes.
 In particular, we prove that any fermion point process
corresponding to a Hermitian kernel may be derived in this way.
}\vspace{3mm}

 \vspace{1.5mm}

 \vspace{1.5mm}


\section{Introduction} Let $X$ be a   locally compact, second countable Hausdorff
topological space. We denote by $\Gamma_X$ the space of locally finite sets (configurations) in $X$.
A point process in $X$ is a probability measure on $\Gamma_X$. From the point of view of classical statistical
mechanics, point processes describe infinite  (interacting) particle systems in continuum.

In the study of point processes, their correlation measures play a crucial role.  Denote by $\Gamma_{X,0}$ the space
of all finite subsets of $X$. One says that a measure $\rho$ on $\Gamma_{X,0}$  is the correlation measure
of a point process $\mu$ if, for each measurable function $G:\Gamma_{X,0}\to[0,+\infty]$, we have (see e.g.\ \cite{Kuna1,Len2,Len3}):
\begin{equation}\label{yufrtyur} \int_{\Gamma_{X,0}}G(\eta) \,\rho(d\eta)=\int_{\Gamma_X}(\mathcal KG)(\gamma)\,\mu(d\gamma),\end{equation}
where the operator $\mathcal K$ is given by
\begin{equation}\label{uyiguyt} (\mathcal KG)(\gamma):=\sum_{\eta\Subset\gamma} G(\eta)\end{equation}
($\eta\Subset\gamma$ denoting that $\eta$ is a finite subset of $\gamma$). It was shown by Lenard \cite{Len1}
that, under a very mild assumption on the correlation measure, it uniquely characterizes a point process.
Furthermore, Lenard \cite{Len3} and Macchi \cite{Macchi} proposed conditions
that are sufficient for a given measure $\rho$ on $\Gamma_{X,0}$ to be the correlation
measure of a point process. Both Lenard and Macchi essentially demanded the local densities derived from
the measure $\rho$ to be non-negative.

Kondratiev and Kuna \cite{Kuna1} treated the $\mathcal K$-transform as an analog of the Fourier transform over the configuration space (see also \cite{Kuna3}). In the case where $X$ is a smooth Riemannian manifold, they defined a $\star$-convolution of functions
on $\Gamma_{X,0}$, so that \begin{equation}\label{hfuytfyut}
\mathcal K(G_1\star G_2)=\mathcal KG_1\cdot\mathcal KG_2,\end{equation} introduced the notion of $\star$-positive definiteness, and proved an analog of the Bochner theorem for point processes.
A spectral approach to this construction, together with a refinement of the local bound satisfied by a measure $\rho$, was proposed in \cite{BKKL}. It should be noted that, in both papers \cite{Kuna1} and \cite{BKKL}, the assumption that $X$ is a smooth manifold was crucial, due to  the use of the Minlos theorem in \cite{Kuna1}, and the  projection
spectral theorem in \cite{BKKL}.

In the first part of this paper (Section \ref{ujftyrr}), we consider a family of commuting Hermitian operators $\mathcal A=(a(\Delta))_{\Delta\in\mathcal B_0(X)} $ indexed
by all measurable, relatively compact sets $\Delta$ in $X$. Such a family of operators may be treated as a quantum stochastic process over the space $X$. We define a class $\mathcal S$ of  ``simple'' functions on $\Gamma_{X,0}$ and,
having fixed the family  $\mathcal A $, introduce corresponding operators
$(\mathcal Q(G))_{G\in\mathcal S}$ such that $\mathcal Q(G_1\star G_2)=\mathcal Q(G_1)\mathcal Q(G_2)$.
We then fix a vector $\Omega$ and say that the family $\mathcal A$ possesses a correlation measure $\rho$ if
$$ (\mathcal Q(G)\Omega,\Omega)=\int_{\Gamma_{X,0}}G(\eta)\,\rho(d\eta).$$
 We prove that, if the family $\mathcal A$  possesses a correlation measure $\rho$ that satisfies some condition of growth, then there exists a point process $\mu $ which has correlation measure $\rho$. Furthermore, the operators
$a(\Delta)$ can be realized as multiplication operators in  $L^2(\Gamma_X,\mu)$.
Thus, $\mu$ can be thought of as the spectral measure of the family $\mathcal A$.
When proving this result, we, in particular,  extend the criterion  of existence of a point process, proved in \cite{Kuna1,BKKL}  to the case
of the topological  space $X$, which is a standard underlying space in the theory of point processes, see e.g.\ \cite{Kal}.

Another tremendous feature of the correlation measure is that it is deeply connected with the normal ordering as it is
known in the quantum field theory. Let us heuristically explain this. Let $\Psi^*(x),\Psi(x)$, $x\in X$, be a representation
of either canonical anticommutation relations (CAR), describing fermions, or canonical  commutation relations
(CAR), describing boson. Then, the operators $a(x):=\Psi^*(x)\Psi(x)$, $x\in X$, describe the particle density and commute,
see e.g. \cite{Goldin1}.  Setting, for each $\Delta\in\mathcal B_0(X)$, $a(\Delta):=\int_\Delta a(x)\,\sigma(dx)$
($\sigma$ being a Radon, non-atomic measure on $X$), we get a family of commuting Hermitian operators.
Let $G^{(n)}$ be a function from $\mathcal S$  such that $G^{(n)}(\eta)=0$ if the number of points in the configuration
$\eta$ is not $n$. Then, one has:
$$ \mathcal Q(G^{(n)})=\frac1{n!}\int_{X^n}\sigma(dx_1)\dotsm\sigma(dx_n)\,G^{(n)}(\{x_1,\dots,x_n\})\Psi^*(x_n)\dotsm
\Psi^*(x_1)\Psi(x_1)\dotsm\Psi(x_n)$$
(normal ordering), so that
\begin{equation}\label{igiuygtu}\frac{n!\, d\rho(\{x_1,\dots,x_n\})}{\sigma(dx_1)\dotsm\sigma(dx_n)}=(\Psi^*(x_n)\dotsm
\Psi^*(x_1)\Psi(x_1)\dotsm\Psi(x_n)\Omega,\Omega),\end{equation}
the expression on the left hand side of \eqref{igiuygtu} being called the $n$-th correlation function.
To the best of our knowledge, heuristic arguments of such type were first given by Menikoff in \cite{M1}, see also
\cite{M2}.

So, in the second part of this paper (Sections \ref{iusfdgh}, \ref{yuyuu}),  we mathematically realize this idea in the case of fermion (determinantal), boson, fermion-like, and boson-like point processes. While fermion and boson point processes have been known since about 1973-1975, when they were introduced by Girard \cite{Girard}, Menikoff \cite{M2}, and Macchi \cite{Macchi} (see also \cite{GY,Sosh,TI,ST} and the references therein), the fermion-like and boson-like point processes first appeared in 2003 in
Shirai and Takahashi's paper \cite{ST}. We also refer to the recent paper \cite{TI}, where, in particular, the case of para-bosons and para-fermions of order 2 is discussed from the  quantum mechanical point of view.

In Section \ref{iusfdgh}, we start with a quasi-free representation of the CAR (CCR, respectively), see e.g.\ \cite{Araki,ArWoods,DA}. Such a representation is completely characterized by a linear, bounded, Hermitian  operator $K$
in $L^2(X,\sigma)$ which satisfies $\pmb0\le K\le\pmb 1$ in the fermion case,  and  $K\ge\pmb0$ in the boson case.
In the case where $X=\mathbb R^d$ and $K$ is a convolution operator, it has been
already shown in \cite{L} that the corresponding particle density has a fermion (boson, respectively) point process as its
spectral measure. We also refer to \cite{BSW}, where a theory of quantum stochastic integration in quasi-free representations of the CAR and CCR was developed   (see also the references therein). 

In this paper, we treat the most general case of the space $X$ and the operator $K$, the only additional
assumption on $K$ being that $K$ is  locally of trace class. The main mathematical (as well as physical) challenge here is to
show that all heuristic arguments coming from physics indeed have a precise mathematical meaning. This is why, at many  steps, we first perform our computations at a heuristic level, and then discuss the mathematical meaning of this procedure.
We observe that  $K$ automatically appears to be an integral operator, and furthermore, with our approach, we do not even have to additionally discuss the problem of the choice of a version of the kernel
$k(x,y)$ of the operator $K$, compare with  \cite[Lemma~1]{Sosh} and \cite[Lemma A4]{GY}.
Thus, we, in particular, show that any fermion process corresponding to a Hermitian operator $K$ can be thought of as the spectral measure
of the family of operators which represent the particle density of a quasi-free representation of the CAR.
Though all our results hold for the complex space $L^2(X\to\mathbb C,\sigma)$, for simplicity of presentation
we only deal with the case of the real space $L^2(X,\sigma)$.

Finally, in Section \ref{yuyuu}, we briefly discuss the family of operators corresponding to fermion-like and (some) boson-like point processes. For a fixed $l\ge2$, we consider a representation of the CAR (CCR, respectively) which is equivalent to the
standard quasi-free representation, but which is based on the orthogonal sum of $2l$ identical copies of the space $L^2(X,\sigma)$ (the standard quasi-free representation being using two copies of this space). The corresponding operators
$\Psi(x)$, $x\in X$,  have the form $\Psi(x)=\sum_{i=1}^l\Psi_i(x)$,  and the particle density is $a(x)=\sum_{i,j=1}^l\Psi_i^*(x)
\Psi_j(x)$. These operators evidently lead to a fermion (boson, respectively) point process. However, we can reduce the
particle density by taking only the ``diagonal elements'' of the double sum: $a^{(l)}(x):=\sum_{i=1}^l\Psi_i^*(x)\Psi_i(x)$.
These operators, in turn, lead to a family of commuting, Hermitian operators $(a^{(l)}(\Delta))_{\Delta\in\mathcal B_0(X)}$,
whose spectral measure is  from the class of point processes discussed in \cite{ST}, and corresponds
to the index $\alpha=-1/l$ ($\alpha=1/l$, respectively) from that paper. Recall that the case $l= 2$ corresponds to the para-fermions (para-bosons, respectively) of order 2, see \cite{TI}. A physical meaning of this procedure of reduction
of the particle density still needs to be clarified.

\section{Correlation measure}\label{ujftyrr}
Let $X$ be a locally compact, second countable Hausdorff
topological space. Recall that such a space is known to be Polish. We denote by $\mathcal{B}(X)$ the Borel $\sigma$-algebra
in $X$, and by $\mathcal{B}_{0}(X)$ the collection of all sets from
$\mathcal{B}(X)$ which are relatively compact.

We define the space of finite multiple configurations in $X$ as
follows:
$$\ddot{\Gamma}_{X,0}:=\bigsqcup_{n\in\mathbb{N}_{0}}\ddot{\Gamma}_{X}^{(n)}.$$
Here, $\mathbb{N}_{0}:=\{0,1,2,\dots\}$, $\ddot{\Gamma}_{X}^{(0)}=\{\varnothing\}$, and for $n\in\mathbb{N}$,
$\ddot{\Gamma}_{X}^{(n)}$ is the factor-space $X^{n}/S_{n}$, where
$S_{n}$ is the group of all permutations of $\{1,\dots,n\}$, which
naturally acts on $X^{n}$:
$$\xi(x_{1},\dots,x_{n})=(x_{\xi(1)},\dots,x_{\xi(n)}),\quad \xi\in
S_{n}.$$ We denote by $[x_{1},\dots,x_{n}]$ the equivalence class
in $\ddot{\Gamma}_{X}^{(n)}$ corresponding to
$(x_{1},\dots,x_{n})\in X^{n}$.

Let $\mathcal{B}(\ddot{\Gamma}_{X}^{(n)})$ denote the image of the
Borel $\sigma$-algebra  $\mathcal{B}(X^{n})$  under the
mapping
$$X^{n}\ni(x_{1},\dots,x_{n})\mapsto[x_{1},\dots,x_{n}]\in\ddot{\Gamma}_{X}^{(n)}.$$
Then, the real-valued measurable functions on
$\ddot{\Gamma}_{X}^{(n)}$ may be identified with the real-valued
$\mathcal{B}_{\mathrm{sym}}(X^{n})$-measurable functions on
$X^{n}$. Here, $\mathcal{B}_{\mathrm{sym}}(X^{n})$ denotes the
$\sigma$-algebra of all sets in $\mathcal{B}(X^{n})$ which are
symmetric, i.e., invariant under the action of $S_{n}$.

For measurable functions
$f_{1},\dots,f_{n}:X\rightarrow\mathbb{R}$, we denote by
$f_{1}\odot\dots\odot f_{n}$ the symmetric tensor product of
$f_{1},\dots,f_{n}$. Since $f_{1}\odot\dots\odot f_{n}$ is
$\mathcal{B}_{\mathrm{sym}}(X^{n})$-measurable, we may consider
$f_{1}\odot\dots\odot f_{n}$ as a measurable function on
$\ddot{\Gamma}_{X}^{(n)}$.

For a function $G:\ddot{\Gamma}_{X,0}\rightarrow \mathbb{R}$, we
denote by $G^{(n)}$ the restriction of $G$ to
$\ddot{\Gamma}_{X}^{(n)}$. Let $\mathcal{S}$ denote the set of all
real-valued functions on $\ddot{\Gamma}_{X,0}$ which satisfy the
following condition: for each $G\in\mathcal{S}$, there is an
$N\in\mathbb{N}$ such that $G^{(n)}=0$ for all $n>N$ and for each
$n\in\{1,\dots,N\}$, $G^{(n)}$ is a finite linear combination of the
functions of the form
$\chi_{\Delta_{1}}\odot\dots\odot\chi_{\Delta_{n}}$, where
$\Delta_{1},\dots,\Delta_{n}\in\mathcal{B}_{0}(X)$ and $\chi_{A}$
denotes the indicator of a set $A$. Note that, by the polarization
identity (e.g.\  \cite[Chapter~2, formula~(2.7)]{BK}), in the above definition it suffices to take
functions of the form $\chi_{\Delta}^{\otimes n}$, where
$\Delta\in\mathcal{B}_{0}(X)$.

Next, we can identify any $[x_{1},\dots,x_{n}]\in\ddot\Gamma_{X,0}$
with the measure $\varepsilon_{x_{1}}+\dots+\varepsilon_{x_{n}}$.
Here, for any $x\in X$, $\varepsilon_{x}$ denotes the Dirac
measure with mass at $x$. We also identify $\{\varnothing\}$ with zero
measure. Through this identification, $\ddot{\Gamma}_{X,0}$
becomes the set of all finite measures on $X$ taking values in
$\mathbb{N}_{0}$.

Next, we introduce the space of finite configurations in $X$,
denoted by $\Gamma_{X,0}$. By definition, $\Gamma_{X,0}$ is the
subset of $\ddot{\Gamma}_{X,0}$ given by
$$\Gamma_{X,0}:=\bigsqcup_{n\in\mathbb{N}_{0}}\Gamma_{X}^{(n)},$$
where $\Gamma_{X}^{(0)}:=\ddot\Gamma_{X}^{(0)}$ and for
$n\in\mathbb{N}$, $\Gamma_{X}^{(n)}$ consists of all
$[x_{1},\dots,x_{n}]\in\ddot{\Gamma}_{X}^{(n)}$ such that
$x_{1},\dots,x_{n}$ are different points in $X$. Hence, each
$[x_{1},\dots,x_{n}]\in\Gamma_{X}^{(n)}$ may be identified either
with the set $\{x_{1},\dots,x_{n}\}$, or with the measure
$\eta=\sum_{i=1}^{n}\varepsilon_{x_{i}}$ satisfying
$\eta(\{x\})\leq 1$ for each $x\in X$. Evidently,
$\Gamma_{X}^{(n)}\in\mathcal{B}(\ddot{\Gamma}_{X}^{(n)})$ and we
denote by $\mathcal{B}(\Gamma_{X}^{(n)})$ the trace
$\sigma$-algebra of $\mathcal{B}(\ddot{\Gamma}_{X}^{(n)})$ on
$\Gamma_{X}^{(n)}$. We also introduce the $\sigma$-algebra
$\mathcal{B}(\Gamma_{X,0})$ on $\Gamma_{X,0}$, whose restriction
to  $\Gamma_{X}^{(n)}$ is $\mathcal{B}(\Gamma_{X}^{(n)})$ for  each
$n\in\mathbb{N}$, and $\{\varnothing\}\in\mathcal{B}(\Gamma_{X,0})$.

Let $\Gamma_{X}$ denote the configuration space over $X$:
$$
\Gamma_{X}:=\{\gamma\subset X:|\gamma\cap\Delta|<\infty
\text{ for each }\Delta\in\mathcal{B}_{0}(X)\}.$$ Here,
$|A|$ denotes the cardinality of a set $A$. Note the evident
inclusion $\Gamma_{X,0}\subset\Gamma_{X}$. We identify each
$\gamma\in\Gamma_{X}$ with the Radon measure
$\gamma=\sum_{x\in\gamma}\varepsilon_{x}$. We introduce the vague
topology on $\Gamma_{X}$ and denote by $\mathcal{B}(\Gamma_{X})$
the Borel $\sigma$-algebra on $\Gamma_{X}$. Note that
$\mathcal{B}(\Gamma_{X})$ is the minimal $\sigma$-algebra on
$\Gamma_{X}$ making all maps
$$\Gamma_{X}\ni\gamma\mapsto\gamma(\Delta)=|\gamma\cap\Delta|\in\mathbb{R},
\quad \Delta\in\mathcal{B}_{0}(X),$$ measurable (cf.\ e.g.\ \cite{Kal}).
Furthermore, the trace $\sigma$-algebra of
$\mathcal{B}(\Gamma_{X})$ on $\Gamma_{X,0}$ coincides with
$\mathcal{B}(\Gamma_{X,0})$.

 For each function $G\in\mathcal S$, we define a measurable function $\mathcal KG:\Gamma_X\to\mathbb R$ by
 \eqref{uyiguyt}.
  Let $\mu$ be a probability measure on $(\Gamma_X,\mathcal B(\Gamma_X))$.
  We say that a measure $\rho$ on $(\Gamma_{X,0},\mathcal B(\Gamma_{X,0}))$ is the correlation measure
 of $\mu$ if  each $G\in\mathcal S$ is integrable with respect to $\rho$, and equality \eqref{yufrtyur} holds for all $G\in\mathcal S$.

We now define a convolution $\star$ as the mapping
$\star:\mathcal{S}\times\mathcal{S}\to\mathcal{S}$ defined
by
\begin{equation}\label{1}
(G_{1}\star G_{2})(\eta):=\sum
G_{1}(\xi_{1}+\xi_{2})G_{2}(\xi_{2}+\xi_{3}),\quad\eta\in \ddot{\Gamma}_{X,0},
\end{equation}
where the summation is over all
$\xi_{1},\xi_{2},\xi_{3}\in\ddot{\Gamma}_{X,0}$ such that
$\xi_{1}+\xi_{2}+\xi_{3}=\eta$. It is easy to see that the
convolution $\star$ is commutative and associative.
Furthermore, we easily see that equality \eqref{hfuytfyut} holds for all $G_1,G_2\in\mathcal S$.

Next, let $F$ be either a real, or  complex Hilbert space and let $D$ be a linear subset of
$F$. Let
$\big(a(\Delta)\big)_{\Delta\in\mathcal{B}_{0}(X)}$ be a
family of Hermitian operators in $F$ such that:
\begin{itemize}
    \item    for each $\Delta\in\mathcal{B}_{0}(X)$, $\operatorname{Dom}\big(a(\Delta)\big)=D$ and
    $a(\Delta)$ maps $D$ into itself;

    \item for any
    $\Delta_1,\Delta_2\in\mathcal{B}_{0}(X)$,   $a(\Delta_1)a(\Delta_2)=
    a(\Delta_2)a(\Delta_1)$;

    \item for any mutually disjoint     $\Delta_{1},\Delta_{2}\in\mathcal{B}_{0}(X)$,
    we have: $a(\Delta_{1}\cup \Delta_{2})=
    a(\Delta_{1})+a(\Delta_{2})$.
\end{itemize}

The family $\big(a(\Delta)\big)_{\Delta\in\mathcal{B}_{0}(X)}$ can
be thought of as a (commutative) quantum stochastic process over a general topological space $X$.

We recursively define the following operators:
\begin{align}
&\mathcal{Q}(\chi_{\Delta_{1}}\odot\dots\odot\chi_{\Delta_{n+1}})
=\frac{1}{(n+1)^{2}}\bigg[\sum_{i=1}^{n+1}a(\Delta_{i})\mathcal{Q}(\chi_{\Delta_{1}}
\odot\dots\odot\check\chi_{\Delta_{i}}\odot\dots\odot
\chi_{\Delta_{n+1}}) \notag\\
&\qquad
-\sum_{i=1}^{n+1}\sum_{j=1,\dots,n+1,\, j\neq
i}\mathcal{Q}\big((\chi_{\Delta_{i}\cap \Delta_{j}})
\odot\chi_{\Delta_{1}}\odot\dots\odot\check{\chi}_{\Delta_{i}}
\odot\dots
\odot\check{\chi}_{\Delta_{j}}\odot\dots\odot\chi_{\Delta_{n+1}}\big)
\bigg],\notag\\
&\quad\quad\Delta_{1},\dots,\Delta_{n+1}\in\mathcal{B}_{0}(X),\ n\in\mathbb{N},\notag\\
&\quad\quad\mathcal{Q}(\chi_{\Delta})=a(\Delta),\quad \Delta\in\mathcal{B}_{0}(X),\label{uiifdyt}
\end{align}
where $\check{\chi}_{\Delta}$ denotes the absence of
$\chi_{\Delta}$. Denote by $\Xi$ the function on $\ddot\Gamma_{X,0}$ given by
$\Xi^{(0)}:=1$, $ \Xi^{(n)}:=0$, $n\in\mathbb{N}$. Let
\begin{equation}\label{ass}\mathcal{Q}(\Xi):= \pmb 1.\end{equation} We
then uniquely define $\mathcal{Q}(G)$ for each $G\in\mathcal{S}$, so
that
$$
\mathcal{Q}(a_{1}G_{1}+a_{1}G_{2})=a_{1}\mathcal{Q}(G_{1})+a_{2}\mathcal{Q}(G_{2}),\quad a_{1},a_{2}\in\mathbb{R},\ G_{1},G_{2}\in\mathcal{S}.
$$
It is not hard to see (e.g., by induction) that, for any
$G_{1},G_{2}\in\mathcal{S}$,
\begin{equation}\label{huy}
\mathcal{Q}(G_{1})\mathcal{Q}(G_{2})=\mathcal{Q}(G_{1}\star
G_{2})
\end{equation}
(compare with \eqref{hfuytfyut}).

We fix any $\Omega\in D$ with $\|\Omega\|_{F}=1$. 
Now, by analogy with \eqref{yufrtyur}, we say that a measure $\rho$ on $\big(\Gamma_{X,0},\mathcal B_0(\Gamma_{X,0})\big)$
is the correlation measure of the family of commuting Hermitian operators $\big(a(\Delta)\big)_{\Delta\in\mathcal{B}_{0}(X)}$ (with respect to the vector $\Omega$) if, for all $G\in\mathcal{S}$,
\begin{equation}\label{14}
(\mathcal Q(G)\Omega,\Omega)_{F}=\int_{\Gamma_{X,0}}G(\eta)\,\rho(d\eta).
\end{equation}

Now, we assume that $\rho$ additionally satisfies:
\begin{description}
        \item[(LB)] {\it Local bound}\,: for each
    $\Delta\in\mathcal{B}_{0}(X)$, there exists $C_{\Delta}>0$ such
    that $$\rho(\Gamma_{\Delta}^{(n)})\leqslant C_{\Delta}^{n},\quad
    n\in\mathbb{N},$$ where
    $\Gamma_{\Delta}^{(n)}:=\{\eta\in\Gamma_{X}^{(n)}\mid\eta\subset\Delta\}$.
    Furthermore, for any sequence
    $\{\Delta_{n}\}_{n\in\mathbb{N}}\in\mathcal{B}_{0}(X)$ such
    that $\Delta_{n}\downarrow\varnothing$, we have $C_{\Delta_{n}}\rightarrow
    0$ as $n\rightarrow\infty$.

\end{description}

\begin{theorem}\label{uytwa} Assume that
a family $\big(a(\Delta)\big)_{\Delta\in\mathcal{B}_{0}(X)}$ of commuting Hermitian operators has
a correlation measure which satisfies  \rom{(LB)}.
For each $G\in\mathcal S$ denote $Q(G):=\mathcal Q(G)\Omega$, and let $\mathfrak{F}$ denote the real Hilbert space obtained as the
closure of the set $\mathfrak{S}:=\{Q(G)\mid G\in\mathcal{S}\}$ in $F$.
For each $\Delta\in\mathcal{B}_{0}(X)$, consider $a(\Delta)$
as an operator in $\mathfrak{F}$ with domain $\mathfrak{S}$. Then,
the operators $a(\Delta)$
are essentially self-adjoint and their closures,
$\tilde{a}(\Delta)$, commute in the sense of their
resolutions of the identity. Furthermore, there exists a unique
probability measure $\mu$ on
$\big(\Gamma_{X},\mathcal{B}(\Gamma_{X})\big)$ whose correlation
measure is $\rho$, the mapping
$$\mathfrak{S}\ni Q(G)\mapsto\big(\mathcal{I}Q(G)\big)(\gamma)
:=\sum_{\eta\Subset\gamma}G(\eta)\in L^{2}(\Gamma,\mu)$$ is
well-defined and extends to a unitary operator
$\mathcal{I}:\mathfrak{F}\rightarrow L^{2}(\Gamma,\mu)$ such that,
under $\mathcal{I}$, $\tilde{a}(\Delta)$ goes over into
the operator of multiplication by $\gamma(\Delta)$.
\end{theorem}

\noindent{\it Proof.} By \eqref{huy} and \eqref{14},
$$\mathcal{S}\times\mathcal{S}\ni(G_{1},G_{2})\mapsto
b_{\rho}(G_{1},G_{2}):=\int_{\Gamma_{X,0}}(G_{1}\star
G_{2})(\eta)\,\rho(d\eta)$$ is a bilinear, positive form on
$\mathcal{S}$. Denote by $\widehat{\mathcal{S}}$ the factorization of
$\mathcal{S}$ consisting of factor-classes
$$\widehat{G}=\{G^{\prime}\in\mathcal{S}:b_{\rho}(G-G^{\prime},G-G^{\prime})=0\},\quad
G\in\mathcal{S}.$$ Define a Hilbert space $\mathcal{H}_{\rho}$ as
the closure of $\widehat{\mathcal{S}}$ in the norm generated by the
scalar product $$(\widehat G_{1},
\widehat G_{2})_{\mathcal{H}_{\rho}}:=b_{\rho}(G_{1},G_{2}).$$

For each $\Delta\in\mathcal{B}_{0}(X)$, we define an operator
$A_{\Delta}$ in $\mathcal{H}_{\rho}$ with domain
$\mathrm{Dom}(A_{\Delta})=\widehat {\mathcal{S}}$ by
\begin{equation}\label{2}
A_{\Delta}\widehat{G}:=\widehat{\chi_{\Delta}\star G},\quad
G\in\mathcal{S}.
\end{equation}
Since $$b_{\rho}(\chi_{\Delta}\star G_{1},G_{2})=b_{\rho}(G_{1},
\chi_{\Delta}\star G_{2}),\quad G_{1},G_{2}\in\mathcal{S},$$ by
\cite[Chapter~5, Section~5, subsec.2]{BK}, the definition (\ref{2}) is indeed
correct, and the operators $A_\Delta$ are symmetric.

Analogously to \cite[Lemma 2]{BKKL}, we easily get the following

\begin{lemma}\label{ifutf}
Each $\widehat{G}\in\widehat{\mathcal{S}}$ is an analytic vector for any
$A_{\Delta}$, $\Delta\in\mathcal{B}_{0}(X)$.
\end{lemma}

By Lemma \ref{ifutf}, for each $\Delta\in\mathcal{B}_{0}(X)$, the closure of
$A_{\Delta}$, denoted by $\tilde A_{\Delta}$, is a self-adjoint
operator in $\mathcal{H}_{\rho}$. Denote by $E_{\Delta}$ its
resolution of the identity. Then, also by Lemma~\ref{ifutf}, for any
$\Delta_{1},\Delta_{2}\in\mathcal{B}_{0}(X)$, the resolutions of
the identity $E_{\Delta_{1}}$ and $E_{\Delta_{2}}$ commute, see e.g.\ \cite[Chapter~5, Theorem~1.15]{BK}.
Therefore, for any
$\Delta_{1},\dots,\Delta_{n}\in\mathcal{B}_{0}(X)$, we can
construct the joint resolution of the identity
$$E_{\Delta_{1},\dots,\Delta_{n}}:=E_{\Delta_{1}}\times\dots \times E_{\Delta_{n}}$$
(see e.g.\ \cite[Chapter~3, Section~1]{BK} for details).

Recall the definition of the function $\Xi$ on $\ddot\Gamma_{X,0}$. Then
$$\nu_{\Delta_{1},\dots,\Delta_{n}}(\cdot):=(E_{\Delta_{1},\dots,
\Delta_{n}}(\cdot)\widehat{\Xi},\widehat{\Xi})_{\mathcal{H}_{\rho}}$$ is a
probability measure on
$(\mathbb{R}^n,\mathcal{B}(\mathbb{R}^{n}))$. Furthermore,
it is clear that
\begin{equation}\label{3}
\{\nu_{\Delta_{1},\dots,\Delta_{n}} \mid
\Delta_{1},\dots,\Delta_{n}\in\mathcal{B}_{0}(X),\ n\in\mathbb{N}\}
\end{equation}
is a consistent family of probability measures. 

 For any $\Delta\in\mathcal{B}_{0}(X)$, denote
$$\Gamma_{\Delta}:=\{\eta\in\Gamma_{X,0}\mid\eta\subset\Delta\},$$
and let $\mathcal{B}(\Gamma_{\Delta})$ be the trace
$\sigma$-algebra of $\mathcal{B}(\Gamma_{X,0})$ on
$\Gamma_{\Delta}$. We introduce a mapping $\mathcal{K}_{\Delta}$, which transforms the
set of all (complex-valued)
functions on $\Gamma_{\Delta}$ into itself, as follows:
\begin{equation}
(\mathcal{K}_{\Delta}G)(\eta):=\sum_{\xi\subset\eta}G(\xi),\quad\eta\in\Gamma_{\Delta}.
\end{equation}
We evidently have:
\begin{equation}\label{6}
\big(\mathcal{K}_{\Delta}(G_{1}\star
G_{2})\big)(\eta)=(\mathcal{K}_{\Delta}G_{1})(\eta)(\mathcal{K}_{\Delta}G_{2})(\eta)
\end{equation}
($G_{1}\star G_{2}$ being given by (\ref{1})). The inverse of
$\mathcal{K}_{\Delta}$ is then given by
\begin{equation}\label{oiyuu}(\mathcal{K}^{-1}_{\Delta}G)(\eta)=\sum_{\xi\subset\eta}
(-1)^{|\eta\backslash\xi|}G(\xi),\quad\eta\in\Gamma_{\Delta}.
\end{equation}

 For any  function $f:\Delta \to\mathbb C$, we define a function $\operatorname{Exp}_\Delta(f,\cdot):\Gamma_{\Delta}\to\mathbb C$ by
 \begin{align*}
 \operatorname{Exp}_\Delta(f,\varnothing):=&1,\\
  \operatorname{Exp}_\Delta(f,\{x_1,\dots,x_n\}):=&f(x_1)\dotsm f(x_n),\quad \{x_1,\dots,x_n\} \in\Gamma_\Delta,\ n\in\mathbb N.
 \end{align*}
By \eqref{oiyuu}, for any $\varphi:\Delta\to\mathbb C$, we have:
\begin{equation}\label{hgufg}
\big(\mathcal K_\Delta^{-1}\exp[\langle\varphi,\cdot\rangle]\big)(\eta)
=\operatorname{Exp}_\Delta(e^\varphi-1,\eta),\quad \eta\in\Gamma_\Delta,
\end{equation}
where $\langle\varphi,\eta\rangle:=\sum_{x\in\eta}\varphi(x)$.

Let $\Delta\in\mathcal{B}_{0}(X)$ be such that
\begin{equation}\label{4}
C_{\Delta}\leqslant\frac{1}{12+\delta},\quad \delta>0
\end{equation}
(see (LB)).
Following the idea of \cite{Kuna1},  we define a set function on  $\mathcal{B}(\Gamma_{\Delta})$ by
\begin{equation}\label{guggy}\mu^{\Delta}(A):=\int_{\Gamma_{\Delta}}(\mathcal{K}_{\Delta}^{-1}
\chi_{A})(\eta)\rho(d\eta),\quad
A\in\mathcal{B}(\Gamma_{\Delta}).\end{equation}
 Since \begin{equation}\label{jhvgj}\sum_{\xi\subset\eta}
1=2^{n}\quad\text{if } |\eta|=n,\end{equation}
 (LB) and (\ref{4}) imply that
$\mu^{\Delta}$ is a signed measure of finite variation.

Let $\Delta_{1},\dots,\Delta_{n}\in\mathcal{B}_{0}(X)$ be  subsets of $\Delta$,  $n\in\mathbb N$, and for simplicity of notations we assume
that these sets are mutually disjoint. Then, by
\eqref{hgufg}--\eqref{jhvgj}, for any  $(y_1,\dots,y_n)\in\mathbb R^n$,
\begin{align}
L(y_1,\dots,y_n):=&\int_{\Gamma_\Delta}\exp\big[\langle i(y_1\chi_{\Delta_1}+\dots+y_n\chi_{\Delta_n}),\eta\rangle\big]\,\mu^
\Delta(d\eta)\notag\\ =& \int_{\Gamma_\Delta}\operatorname{Exp}_\Delta
\big((e^{iy_1}-1)\chi_{\Delta_1}+\dots+(e^{iy_n}-1)\chi_{\Delta_n},\eta
\big)\rho(d\eta).\label{ufytfy}
\end{align}
Using (LB), \eqref{6}, \eqref{4}, and   \eqref{ufytfy}, we conclude that the function $L:\mathbb R^n\to\mathbb C$ is positive definite in the sense of the Fourier analysis on $\mathbb R^n$. Hence, $L$ is the Fourier transform of a probability measure  on $\mathbb R^n$. Therefore,  under the mapping
$$\Gamma_{\Delta}\ni\eta\mapsto\big(\eta(\Delta_{1}),\dots,\eta(\Delta_{n})\big)\in\mathbb{R}^{n},$$
the image of the signed measure $\mu^{\Delta}$ is a probability
measure on $(\mathbb R^n,\mathcal B(\mathbb R^n))$, which we denote by $\mu^\Delta_{\Delta_1,\dots,\Delta_n}$.

Using \eqref{6}, \eqref{oiyuu}, and \eqref{4}, for any $y^{(1)},\dots,y^{(k)}\in\mathbb R^n$, $k\in\mathbb N$, we have:
\begin{align}
\int_{\mathbb R^n}\prod_{i=1}^k (x,y^{(i)})_{\mathbb R^n}
  \, d\nu_{\Delta_1,\dots,\Delta_n}(x)
&=\bigg(\prod_{i=1}^k\bigg(
\sum_{j=1}^n A_{\Delta_j}y_j^{(i)}\bigg)\hat\Xi,\hat\Xi\bigg)_{\mathcal H_{\rho}}
\notag\\
&=\int_{\Gamma_{X,0}}\bigg(\sum_{j=1}^n y_j^{(1)}\chi_{\Delta_j}\bigg)
\star \dots\star\bigg(\sum_{j=1}^n y_j^{(k)}\chi_{\Delta_j}\bigg)\,\rho(d\eta)\notag\\
& =\int_{\mathbb R^n}\prod_{i=1}^k (x,y^{(i)})_{\mathbb R^n}
  \, d\mu_{\Delta_1,\dots,\Delta_n}^\Delta(x).
\label{erees}
\end{align}
Furthermore,  it follows from  the proof of
Lemma~\ref{ifutf} that  there exists a constant $R>0$ such that
\begin{equation}\bigg| \int_{\Gamma_{X,0}}\bigg(\sum_{j=1}^n y_j^{(1)}\chi_{\Delta_j}\bigg)
\star \dots\star\bigg(\sum_{j=1}^n y_j^{(k)}\chi_{\Delta_j}\bigg)\,\rho(d\eta)\bigg|\le R^n n! \prod_{i=1}^k
\| y^{(i)}\|_{\mathbb R^n}.\label{gig}\end{equation}
Hence, by the theorem on  uniqueness of the solution of a moment problem (e.g.\ \cite[Chapter~5, Theorem~2.1 and Remark~3]{BK}), we conclude from  \eqref{erees} and \eqref{gig} that
$$\nu_{\Delta_1,\dots,\Delta_n}=\mu_{\Delta_1,\dots,\Delta_n}^\Delta.$$

We also observe that the sets
\begin{align}\label{8}
\{\eta\in\Gamma_{X,0}\mid&\big(\eta(\Delta_{1}),\dots,\eta(\Delta_{n})\big)\in
B_{n}\},\nonumber\\
&B_{n}\in\mathcal{B}(\mathbb{R}^{n}),\ \Delta_{1},\dots,\Delta_{n}\in\mathcal{B}_{0}(X),\
 \Delta_{1}\cup\dotsm\cup\Delta_{n}\subset\Delta,\ n\in\mathbb{N},
\end{align}
 generate the $\sigma$-algebra
$\mathcal{B}(\Gamma_{\Delta})$. Hence,  $\mu_{\Delta}$ is
a probability measure on
$\big(\Gamma_{\Delta},\mathcal{B}(\Gamma_{\Delta})\big)$.

Next, let $\Delta^{\prime}\in\mathcal{B}_{0}(X)$ be such that
$\Delta^{\prime}\subset\Delta$. As usual, we identity
$\mathcal{B}(\Gamma_{\Delta^{\prime}})$ with the
sub-$\sigma$-algebra of $\mathcal{B}(\Gamma_{\Delta})$ generated
by the sets of the form (\ref{8}) where
$\Delta_{1},\dots,\Delta_{n}$ are subsets of $\Delta^{\prime}$.
Then it follows from the above that $\mu_{\Delta^{\prime}}$ is the
restriction of $\mu_{\Delta}$ to
$\mathcal{B}(\Gamma_{\Delta^{\prime}})$.

Now, we will show that there exists a random measure $M$ on $X$
such that, for any
$\Delta_{1},\dots,\Delta_{n}\in\mathcal{B}_{0}(X)$, $n\in\mathbb{N}$,
the distribution of $\big(M(\Delta_{1}),\dots,M(\Delta_{n})\big)$
is $\nu_{\Delta_{1},\dots,\Delta_{n}}$ (see e.g.\ \cite{Kal} for details
on random measures).

By (LB), for any $x\in X$, there exists an open neighborhood of
$x$, denoted by $\Delta(x)$, such that
$\Delta(x)\in\mathcal{B}_{0}(X)$ and $C_{\Delta(x)}\leqslant 1/(12+\delta)$.
Therefore, for any $\Delta\in\mathcal{B}_{0}(X)$, there exist
mutually disjoint sets
$\Delta_{1},\dots,\Delta_{m}\in\mathcal{B}_{0}(X)$, $m\in\mathbb{N}$,
such that
$\Delta=\Delta_{1}\cup\dots\cup\Delta_{m}$, $C_{\Delta_{i}}\leqslant
1/(12+\delta)$,  $i=1,\dots,m$. By the proved above,
$\nu_{\Delta_{i}}\big([0,+\infty)\big)=1$, $i=1,\dots,m$. Hence,
$\nu_{\Delta_{1},\dots,\Delta_{m}}([0,+\infty)^{m})=1$. By Lemma~\ref{ifutf},
 for each $A\in\mathcal{B}(\mathbb{R})$,
$$\nu_{\Delta}(A)=\nu_{\Delta_{1}\cup\dots\cup\Delta_{m}}(A)
=\int_{[0,+\infty)^{m}}\chi_{A}(x_{1}+\dots+x_{m})\,d\nu_{\Delta_{1},
\dots,\Delta_{m}}(x_{1},\dots,x_{m}).$$ Hence,
$\nu_{\Delta}\big([0,+\infty)\big)=1$, and so, for any
$\Delta_{1},\dots,\Delta_{n}\in\mathcal{B}_{0}(X)$, $n\in\mathbb{N}$,
$$\nu_{\Delta_{1},\dots,\Delta_{n}}\big([0,+\infty)^{n}\big)=1.$$

Next, it is also clear from Lemma~\ref{ifutf}  that, for any disjoint
$\Delta_{1},\Delta_{2}\in\mathcal{B}_{0}(X)$,
$$\nu_{\Delta_{1},\Delta_{2},\Delta_{1}\cup\Delta_{2}}
\big(\{(x,y,z)\in\mathbb{R}^{3}\mid x+y=z\}\big)=1.$$

Finally, let $\Delta_{n}\in\mathcal{B}_{0}(X),\,\,n\in\mathbb{N}$,
be such that $\Delta_{n}\downarrow\varnothing$. We state that
$\nu_{\Delta_{n}}$ weakly converges to $\varepsilon_{0}$. By (A2),
without loss, we may assume that $C_{\Delta_{1}}\leqslant1/(12+\delta)$.
Then, each $\nu_{\Delta_{n}}$ is concentrated on the set
$\mathbb{N}_{0}$. Hence, it is enough to show that
$\nu_{\Delta_{n}}(\mathbb{N})\rightarrow0$ as
$n\rightarrow\infty$. But this holds since $\nu_{\Delta_{n}}$ is
the distribution of the random variable $\eta(\Delta_{n})$ under
$\mu^{\Delta_1}$.

Now, by \cite[Theorem 5.4]{Kal}, there indeed exists a random measure $M$ on $X$ as
described above. Furthermore, we already know that, for any $x\in X$, there exists
an  open neighborhood of $x$, denoted by $\Delta(x)$, such that
$\Delta(x)\in\mathcal{B}_{0}(X)$ and the restriction of $M$ to $\Delta(x)$ is
concentrated on  $\Gamma_{\Delta(x)}$.  Hence, the random measure $M$ is a.s.\ concentrated on
$\Gamma_{X}$. Letting $\mu$ denote the distribution of $M$ on $\Gamma_{X}$, we
obtain a unique probability measure on
$\big(\Gamma_{X},\mathcal{B}(\Gamma_{X})\big)$ whose ``finite-dimensional
distributions'' are given through the measures (\ref{3}).

Let again $\Delta\in\mathcal{B}_{0}(X)$ be such that (\ref{4}) is satisfied. As
usual, we identify $\mathcal{B}(\Gamma_{\Delta})$ as a
sub-$\sigma$-algebra of $\mathcal{B}(\Gamma_{X,0})$. Then, for any
$G_{1},G_{2}\in\mathcal{S}$ which, restricted to $\Gamma_{X,0}$,
are $\mathcal{B}(\Gamma_{\Delta})$-measurable, we have:
\begin{align}\label{10}
&\int_{\Gamma_{X,0}}(G_{1}\star G_{2})(\eta)\,\rho(d\eta)\nonumber\\
&\qquad=\int_{\Gamma_{\Delta}}(G_{1}\star
G_{2})(\eta)\rho(d\eta)\nonumber\\
&\qquad=\int_{\Gamma_{\Delta}}(\mathcal{K}_{\Delta}G_{1})(\eta)
(\mathcal{K}_{\Delta}G_{2})(\eta)\,\mu^{\Delta}(d\eta)\nonumber\\
&\qquad=\int_{\Gamma_{X}}\bigg(\sum_{\eta\Subset\gamma}G_{1}(\eta)\bigg)
\bigg(\sum_{\eta\Subset\gamma}G_{2}(\eta)\bigg)\,\mu(d\eta).
\end{align}

Next, any $G\in\mathcal{S}$ can be represented as
$G=\sum_{j=1}^{k}G_{j}$, where $k\in\mathbb{N}$, each
$G_{j}$ belongs to $\mathcal{S}$ and, restricted to $\Gamma_{X,0}$, is
$\mathcal{B}(\Gamma_{\Delta_{j}})$-measurable with
$\Delta_{j}\in\mathcal{B}_{0}(X)$, $C_{\Delta_{j}}\leqslant1/(12+\delta)$.
Hence, by (\ref{10}), for any $G_{1},G_{2}\in\mathcal{S}$,
\begin{equation}\label{kjgkgiugy}\int_{\Gamma_{X,0}}(G_{1}\star G_{2})(\eta)\,\rho(d\eta)
=\int_{\Gamma_{X}}\bigg(\sum_{\eta\Subset\gamma}G_{1}(\eta)\bigg)
\bigg(\sum_{\eta\Subset\gamma}G_{2}(\eta)\bigg)\mu(d\eta).\end{equation}
Define the mapping
\begin{equation}\label{9}
\widehat{\mathcal{S}}\ni\widehat{G}\mapsto(\mathcal{K}\widehat{G})(\gamma):=\sum_{\eta\Subset\gamma}G(\eta).
\end{equation}
Then, by \eqref{kjgkgiugy}, $\mathcal K$
extends to an isometry of $\mathcal{H}_{\rho}$ into $L^{2}(\Gamma,\mu)$.
Furthermore, the image of $\mathcal{K}$ is evidently dense in
$L^{2}(\Gamma,\mu)$, and so $\mathcal{K}$ is a unitary operator.

For each $\Delta\in\mathcal{B}_{0}(X)$,
$$\mathcal{K}(\widehat{\chi_{\Delta}\star G})(\gamma)
=\gamma(\Delta)(\mathcal{K}\widehat{G})(\gamma),\quad G\in\mathcal{S},\ \gamma\in\Gamma_{X}.$$
Therefore,
$\tilde{A}_{\Delta}$ goes over, under $\mathcal{K}$, into the operator of
multiplication by $\gamma(\Delta)$.

Finally, we can construct a unitary operator $I:\mathfrak F\to \mathcal H_\rho$ by setting
$IQ(G):=\widehat G$. Then, from the proved above,   we get the conclusion of the theorem, except for the statement about the
uniqueness of a measure $\mu$, whose  correlation measure is $\rho$.
 But the uniqueness of such $\mu$ follows from \cite{Len1}  (in fact,
the uniqueness can also be derived from the above arguments).  \quad $\square$\vspace{2mm}

It is clear that any
 correlation measure $\rho$ satisfies  the following condition:
\begin{description}
    \item[(N)] {\it Normalization}\,: $\rho(\Gamma_{X}^{(0)})=1$.
\end{description}

It follows from \eqref{huy} and \eqref{14} (or \eqref{hfuytfyut}) that any correlation measure  $\rho$ also satisfies:

\begin{description}
  \item[(PD)] {\it $\star$-positive definiteness}\,: For each
    $G\in\mathcal{S}$: $$\int_{\Gamma_{X,0}}(G\star G)(\eta)\,\rho(d\eta)\geqslant
    0.$$
\end{description}

From (the proof of) Theorem~\ref{uytwa}, we easily conclude the following criterion
of existence of a point process, which generalizes  \cite[Theorem~6.5]{Kuna1} and  \cite[Theorem~2]{BKKL}.

\begin{corollary}\label{uitgyuit}
Let $\rho$ be a measure on
$\big(\Gamma_{X,0},\mathcal{B}(\Gamma_{X,0})\big)$ satisfying
\rom{(N)}, \rom{(PD)}, and \rom{(LB)}. Then, there exists a unique probability measure on
$\big(\Gamma_{X},\mathcal{B}(\Gamma_{X})\big)$ which has $\rho$ as
correlation measure.
\end{corollary}

\section{Particle densities in quasi-free representations of the CAR and CCR}\label{iusfdgh}
Let $X$ be a topological space as in Section \ref{ujftyrr}. Let $\sigma$ be a
non-atomic Radon measure on $(X,\mathcal{B}(X))$. We denote by $H$
the real space $L^{2}(X,\sigma)$. For an integral operator $I$ in $H$, we denote by $\mathcal
N(I)$ the kernel of $I$.

Let $K$ be a linear bounded operator in $H$ which satisfies the
following assumptions:
\begin{itemize}
    \item $K$ is symmetric and $\pmb 0 \leqslant K \leqslant \pmb 1$;
    \item for each $\Delta\in \mathcal{B}_{0}(X)$, the operator
    $P_{\Delta}KP_{\Delta}$ is of trace class. Here, $P_{\Delta}$
    denotes the operator of multiplication by $\chi_{\Delta}$.
\end{itemize}

Denote $K_{1}:=\sqrt{K}$. For each $\Delta\in \mathcal{B}_{0}(X)$,
$$P_{\Delta}K_{1}(P_{\Delta}K_{1})^{*}=P_{\Delta}KP_{\Delta}.$$
Therefore, the operator $P_{\Delta}K_{1}$ is of Hilbert--Schmidt class. Hence,
$P_{\Delta}K_{1}$ is an integral operator, whose kernel $\mathcal N(P_\Delta
K_1)$ belongs to $L^{2}(X^{2}, \sigma^{2})$. This implies that $K_{1}$ is an
integral operator, whose kernel satisfies
\begin{equation}\label{*}
\int_{\Delta}\int_{X} \mathcal N(K_1) (x,y)^2\,\sigma(dx)\,\sigma(dy)<
\infty,\quad \Delta\in \mathcal{B}_{0}(X).
\end{equation}
Note also that the kernel $\mathcal N(K_1)$ is symmetric.

Thus, $K$ is an integral operator, whose kernel is given by
$$k(x,y):=\mathcal N(K)(x,y)=\int_{X}\mathcal N(K_1)(x,z)\mathcal N(K_1)(z,y)
\,\sigma(dz).$$ By (\ref{*}),
for any $\Delta\in \mathcal{B}_{0}(X)$, we get:
\begin{align*}
\int_{\Delta}k(x,x)\,\sigma(dx)&=
\int_{\Delta}\int_{X}\mathcal N(K_1)(x,y) \mathcal N(K_1)(y,x)\,\sigma(dy)\,\sigma(dx)\\
&=\int_{\Delta}\int_{X}\mathcal N(K_1)(x,y)^2\,\sigma(dx)\,\sigma(dy)<\infty.
\end{align*}
Note that the kernel  $\mathcal N(K_1)(x,y)$  is defined up to a set of
$\sigma^{\otimes2}$-measure 0 in $X^{2}$, but the value
$\int_{X}k(x,x)\,\sigma(dx)$ is independent of the choice of a version of
$\mathcal N(K_1)$.

Now, for a fixed $x\in X$, we define the function
$\varkappa_{1,x}:X\rightarrow\mathbb{R}$  by $$\varkappa_{1,x}(y):=\mathcal
N(K_1)(x,y),\quad y\in X.$$ By (\ref{*}),
\begin{equation}\label{222}
\varkappa_{1,x} \in L^{2}(X,\sigma) \quad \text{for $\sigma$-a.a.\ }x\in X.
\end{equation}
We also define the linear bounded operator $K_{2}:=(\pmb 1-K)^{1/2}$. Though
$K_{2}$ is not an integral operator, we will heuristically use
$\varkappa_{2,x}$, $x\in X$, to denote the ``function''
$\varkappa_{2,x}(y):=\mathcal N(K_2)(x,y)$, where $\mathcal N(K_2)(x,y)$ is the
``kernel'' of $K_{2}$.

For a real separable Hilbert space $\mathcal H$, we denote by
$\mathcal A\mathcal F(\mathcal H)$ the antisymmetric Fock space over $\mathcal H$:
$$\mathcal A\mathcal F(\mathcal H):=\bigoplus_{n=0}^{\infty}\mathcal A\mathcal F^{(n)}(\mathcal
H).$$  Here,
$\mathcal A\mathcal F^{(0)}(\mathcal H):=\mathbb{R}$ and for $n\in
\mathbb{N}$ $\mathcal A\mathcal F^{(n)}(\mathcal H):=\mathcal
H^{\wedge n}n!\,$, where $\wedge$ stands for the antisymmetric
tensor product and $n!$ is a normalizing factor, so that, for any
$f^{(n)}\in \mathcal A\mathcal F^{(n)}(\mathcal H)$,
$$\|f^{(n)}\|^{2}_{\mathcal A\mathcal F^{(n)}(\mathcal H)}=\|f^{(n)}\|^{2}_{\mathcal H^{\wedge
n}}\,n!\,.$$ We denote by $\mathcal A\mathcal
F_{\mathrm{fin}}(\mathcal H)$ the subset of $\mathcal A\mathcal
F(\mathcal H)$ consisting of all elements
$f=(f^{(n)})_{n=0}^{\infty}\in\mathcal A\mathcal F(\mathcal H)$
for which $f^{(n)}=0$, $n\geqslant N$, for some $N\in \mathbb{N}$. We
endow $\mathcal A\mathcal F_{\mathrm{fin}}(\mathcal H)$ with the
topology of the topological direct sum of the spaces $\mathcal
A\mathcal F^{(n)}(\mathcal H)$. Thus, the convergence in $\mathcal A\mathcal F_{\mathrm{fin}}(\mathcal H)$ means uniform boundedness and
coordinate-wise convergence.

For $g\in\mathcal H$, we denote by $\Phi(g)$ and $\Phi^{*}(g)$ the annihilation
and creation operators in $\mathcal A\mathcal F(\mathcal H)$, respectively.
These are linear continuous operators in $\mathcal A\mathcal
F_{\mathrm{fin}}(\mathcal H)$ defined through the formulas
\begin{align*}
\Phi(g)\,h_{1}\wedge\dots\wedge
h_{n}:=&\sum_{i=1}^{n}(-1)^{i+1}(g,h_{i})_{\mathcal H}\, h_{1}\wedge\dots\wedge
h_{i-1}\wedge \check{h}_{i}\wedge
h_{i+1}\wedge\dots\wedge h_{n},\\
\Phi^{*}(g)\,h_{1}\wedge\dots\wedge h_{n}:=&g\wedge h_{1}\wedge\dots\wedge
h_{n},
\end{align*}
where $h_{1},\dots,h_{n}\in \mathcal H$.

We now set $\mathcal H:=H_{1}\oplus H_{2}$, where $H_{1}$ and
$H_{2}$ are two copies of $H$. For $f\in H$, we denote
$$\Phi_{1}(f):=\Phi(f,0),\quad \Phi_{2}(f):=\Phi(0,f),$$ and analogously
$\Phi_{i}^{*}(f)$, $i=1,2$. We set, for each $f\in H$,
\begin{align}
\Psi(f):=&\Phi_{2}(K_{2}f)+\Phi_1^{*}(K_{1}f),\notag\\
\Psi^{*}(f):=&\Phi_{2}^{*}(K_{2}f)+\Phi_{1}(K_{1}f).\label{vfsiu}
\end{align}
The operators $\{\Psi(f),\Psi^{*}(f)\mid f\in \mathcal H\}$ satisfy the CAR:
\begin{align}
[\Psi(f),\Psi(g)]_+&=[\Psi^*(f),\Psi^*(g)]_+=0,\notag\\
[\Psi^*(f),\Psi(g)]_+&=(f,g)_H\,\pmb 1,\quad f,g\in H,\label{uiyftur}
\end{align}
where $[A,B]_+:=AB+BA$. This representation of the CAR is called
quasi-free. The so-called $n$-point functions  of this representation have the structure
\begin{equation}\label{uytfrtyrfrr}
(\Psi^*(f_n)\dotsm\Psi^*(f_1)\Psi(g_1)\dotsm\Psi(g_m)\Omega,\Omega)
_{\mathcal{AF}(\mathcal
H)}=\delta_{n,m}\operatorname{det}[(Kf_i,g_j)_H]_{i,j=1}^n.
\end{equation}
Here, $\Omega:=(1,0,0,\dots)$ is the vacuum vector in $\mathcal {AF}(\mathcal H)$.

We have the following heuristic representation:
\begin{align*}
\Psi(f)&=\int_{X}\sigma(dx)\,f(x)\Psi(x)\\
&=\int_{X}\sigma(dx)\big((K_{2}f)(x)\Phi_{2}(x)+(K_{1}f)(x)\Phi_{1}^{*}(x)\big)\\
&=\int_{X}\sigma(dx)\bigg(\Phi_{2}(x)\int_{X}\sigma(dy)\, \mathcal
N(K_2)(x,y)f(y)+\Phi_{1}^{*}(x)\int_{X}\sigma(dy)\,\mathcal N(K_1)
(x,y)f(y)\bigg)\\
&=\int_{X}\sigma(dy)\,f(y)\bigg(\int_{X}\sigma(dx)\,\mathcal N(K_2)(x,y)
\Phi_{2}(x)+\int_{X}\sigma(dx)\,\mathcal N(K_1)(x,y)\Phi_{1}^{*}(x)\bigg)\\
&=\int_{X}\sigma(dx)\,f(x)\big(\Phi_{2}(\varkappa_{2,x})+\Phi_{1}^{*}(\varkappa_{1,x})
\big),
\end{align*}
and analogously
\begin{align*}
\Psi^{*}(f)&=\int_{X}\sigma(dx)\,f(x)\Psi^{*}(x)\\
&=\int_{X}\sigma(dx)\,f(x)\big(\Phi_{2}^{*}(\varkappa_{2,x})
+\Phi_{1}(\varkappa_{1,x})\big).
\end{align*}
Hence, for $x\in X$,
\begin{align*}
\Psi(x)&=\Phi_{2}(\varkappa_{2,x})+\Phi_{1}^{*}(\varkappa_{1,x}),\\
\Psi^{*}(x)&=\Phi_{2}^{*}(\varkappa_{2,x})+\Phi_{1}(\varkappa_{1,x}).
\end{align*}

We now heuristically define
\begin{align*}
a(x):&=\Psi^{*}(x)\Psi(x),\quad x\in X,\\
a(\Delta):&=\int_{\Delta}\sigma(dx)\,a(x)\\
&=\int_\Delta\sigma(dx)\, \big(\Phi_{2}^{*}(\varkappa_{2,x})+\Phi_{1}(\varkappa_{1,x})\big)
\big(\Phi_{2}(\varkappa_{2,x})+\Phi_{1}^{*}(\varkappa_{1,x})\big), \quad \Delta\in \mathcal B_0(X).
\end{align*}

We will now show that it is, in fact, possible to realize
$a(\Delta)$, $\Delta\in \mathcal B_0(X)$, as linear continuous operators on
$\mathcal A\mathcal F_{\mathrm{fin}}(\mathcal H)$.
We first look at the operator
$$a^{+}(\Delta):=\int_{\Delta}\sigma(dx)\,\Phi_{2}^{*}(\varkappa_{2,x})\Phi_{1}^{*}(\varkappa_{1,x}).$$
We heuristically have:
\begin{equation}\label{333}
a^{+}(\Delta)h_{1}\wedge\dots\wedge
h_{n}=\int_{\Delta}\sigma(dx)\,\varkappa_{2,x}\wedge\varkappa_{1,x}\wedge
h_{1}\wedge\dots\wedge h_{n}.
\end{equation}
Here and below, we identify $\varkappa_{1,x}$ with
$(\varkappa_{1,x},0)$ and $\varkappa_{2,x}$ with
$(0,\varkappa_{2,x})$.

To make sense out of (\ref{333}), we need to
show that the informal expression
\begin{equation}\label{yhufgu}\int_{\Delta}\sigma(dx)\,\varkappa_{2,x}\otimes\varkappa_{1,x}\end{equation}
identifies an element of the Hilbert space $\mathcal H^{\otimes 2}$.
So, take any $u\in H_2$ and $v\in H_1$. Then
\begin{align}\label{444}
\bigg(\int_{\Delta}\sigma(dx)\,\varkappa_{2,x}\otimes\varkappa_{1,x},
u\otimes
v\bigg)_{\mathcal H^{\otimes 2}}
&=\int_{\Delta}\sigma(dx)\,(\varkappa_{2,x}, u)_{H_2}(\varkappa_{1,x},v_1)_{H}\nonumber\\
&=\int_{\Delta}\sigma(dx)\,\big(K_{2}u\big)(x)\big(K_{1}v\big)(x)\nonumber\\
&=\int_{X}\sigma(dx)\,u(x)\big(K_{2}P_{\Delta}K_{1}v\big)(x).
\end{align}
Since $P_{\Delta}K_{1}$ is of Hilbert--Schmidt
class, hence so is $K_{2}P_{\Delta}K_{1}$. Therefore, the operator
$K_{2}P_{\Delta}K_{1}$ has a kernel $\mathcal N(K_{2}P_{\Delta}K_{1})$, which belongs to $H^{\otimes 2}$. Therefore, we continue (\ref{444}) as follows:
\begin{align}
&=\int_{X}\sigma(dx)\,u(x)\int_{X}\sigma(dy)\,\mathcal N(K_{2}P_{\Delta}K_{1})(x,y)v(y)\nonumber\\&
=(\mathcal N(K_{2}P_{\Delta}K_{1})_{2,1},u\otimes v)_{\mathcal H^{\otimes 2}}.\label{555}
\end{align}
Here, $\mathcal N(K_{2}P_{\Delta}K_{1})_{2,1}$ is the element of the space $\mathcal
H^{\otimes 2}$ which belongs to
its subspace $H_{2}\otimes H_{1}$ and coincides in it with $\mathcal N(K_{2}P_{\Delta}K_{1})$. Let also $\mathcal N(K_{2}P_{\Delta}K_{1})_{2,1}^\wedge$ denote the orthogonal projection of
$\mathcal N(K_{2}P_{\Delta}K_{1})_{2,1}$ onto $\mathcal
H^{\wedge 2}$.

Thus, the rigorous definition of $a^{+}(\Delta)$ is as follows:
\begin{equation}\label{jhvsu}
a^{+}(\Delta)h_{1}\wedge\dots\wedge h_{n}=\mathcal N(K_{2}P_{\Delta}K_{1})_{2,1}^\wedge\wedge
h_{1}\wedge\dots\wedge h_{n},
\end{equation}
i.e., $a^{+}(\Delta)$ is the creation by $\mathcal N(K_{2}P_{\Delta}K_{1})_{2,1}^\wedge$. The
$a^{+}(\Delta)$ is evidently a linear continuous operator on
$\mathcal A\mathcal F_{\mathrm{fin}}(\mathcal H)$.

Next, the operator
$$a^{-}(\Delta):=\int_{\Delta}\sigma(dx)\,\Phi_{1}(\varkappa_{1,x})\Phi_{2}(\varkappa_{2,x})$$
should be rigorously understood as the restriction to $\mathcal
A\mathcal F_{\mathrm{fin}}(\mathcal H)$ of the adjoint operator
$\big(a^{+}(\Delta)\big)^{+}$. Hence
\begin{equation}\label{777}
a^{-}(\Delta)h_{1}\wedge\dots\wedge h_{n}=n(n-1)\big(\mathcal
N(K_{2}P_{\Delta}K_{1})_{2,1},h_{1}\wedge\dots\wedge h_{n}\big)_{\mathcal
H^{\otimes 2}},
\end{equation}
where the scalar product is taken in the first two ``variables''. Therefore,
\begin{align*}
a^{-}(\Delta)h_{1}\wedge\dots\wedge
h_{n}=&\sum_{i,j=1,\dots,n,\,i\neq
j}(-1)^{i+j+\chi_{\{i<j\}}(i,j)}\big(\mathcal N(K_{2}P_{\Delta}K_{1}),h_{i}^{(2)}\otimes
h_{j}^{(1)}\big)_{H^{\otimes 2}}\\
&\times h_{1}\wedge\dots\wedge\check{h}_{i}\wedge\dots\wedge\check{h}_{j}\wedge\dots\wedge
h_{n}.
\end{align*}
Here and below, we use the notation $h_{i}=\big(h_{i}^{(1)},h_{i}^{(2)}\big)$.
Note also that
$$\big(\mathcal N(K_{2}P_{\Delta}K_{1}),h_{i}^{(2)}\otimes h_{j}^{(1)}\big)_
{H^{\otimes 2}} =\big(h_{i}^{(2)},K_{2}P_{\Delta}K_{1}h_{j}^{(1)}\big)_{H}.$$

For the operator
$$a_{1}^{0}(\Delta):=\int_{\Delta}\sigma(dx)\,\Phi_{1}(\varkappa_{1,x})\Phi_{1}^{*}(\varkappa_{1,x}),$$
we have (recall (\ref{222})\big):
\begin{align}
&a_{1}^{0}(\Delta)h_{1}\wedge\dots\wedge
h_{n}\notag\\
&\qquad =\int_{\Delta}\|\varkappa_{1,x}\|_{H}^{2}\,\sigma(dx)\,h_{1}\wedge\dots\wedge
h_{n}\notag\\&\qquad\quad
-\sum_{i=1}^{n}h_{1}\wedge\dots\wedge
h_{i-1}\wedge\bigg(\int_{\Delta}\sigma(dx)\,(\varkappa_{1,x},h_{i})
_{\mathcal H}\varkappa_{1,x}\bigg)
\wedge h_{i+1}\wedge\dots\wedge h_{n}.\label{888}
\end{align}
For any $u,v\in H$, we have:
\begin{align}
\bigg(\int_{\Delta}\sigma(dx)\,(\varkappa_{1,x},u)_{H}\,\varkappa_{1,x},v\bigg)_{H}
&=\int_{\Delta}\sigma(dx)\,\big(K_{1}u\big)(x)\big(K_{1}v\big)(x)\nonumber\\
&=\big(K_{1}P_{\Delta}K_{1}u,v\big)_{H}\,.\label{999}
\end{align}

For any linear operator $A$ on $\mathcal H$, we define the second
quantization of $A$, denoted by $d\Gamma(A)$, as the linear
continuous operator on $\mathcal A\mathcal
F_{\mathrm{fin}}(\mathcal H)$ given by
\begin{align*}
d\Gamma(A)\upharpoonright\mathcal A\mathcal F^{(0)}(\mathcal
H)&=\pmb 0,\\
d\Gamma(A)\upharpoonright\mathcal A\mathcal F^{(n)}(\mathcal
H)&=A\otimes\pmb 1\otimes\dots\otimes\pmb 1+\pmb 1\otimes A\otimes\pmb
1\otimes\dots\otimes\pmb 1\\
&\quad+\dots+\pmb 1\otimes\dots\otimes\pmb 1\otimes A,\quad n\in
\mathbb{N}.
\end{align*}
We identify the operator $K_{1}P_{\Delta}K_{1}$ in $H$ with the
operator $K_{1}P_{\Delta}K_{1}\oplus\pmb0$ in $\mathcal H$. Then, by virtue
of (\ref{888}) and  (\ref{999}), we define:
\begin{equation*}
a_1^{0}(\Delta):=\int_{\Delta}\|\varkappa_{1,x}\|_{H}^{2}\,\sigma(dx)\,\pmb
1-d\Gamma\big(K_{1}P_{\Delta}K_{1}\big).
\end{equation*}

Analogously, we conclude that the operator $a_{2}^{0}(\Delta)$,
which is heuristically given by
$$a_{2}^{0}(\Delta)=\int_{\Delta}\sigma(dx)\,\Phi_{2}^{*}(\varkappa_{2,x})\Phi_{2}(\varkappa_{2,x}),$$
can be rigorously defined as
$$a_{2}^{0}(\Delta):=d\Gamma\big(K_{2}P_{\Delta}K_{2}\big),$$ where we
identified the operator $K_{2}P_{\Delta}K_{2}$ in $H$ with the
operator $\pmb 0\oplus K_{2}P_{\Delta}K_{2}$ in $\mathcal H$.

We next define, for $\Delta\in \mathcal B_0(X)$,
\begin{align}
a^0(\Delta):=&a_{1}^{0}(\Delta)+a_{2}^{0}(\Delta)\notag\\
=&\int_{\Delta}\|\varkappa_{1,x}\|_{H}^{2}\,\sigma(dx)\,\pmb 1
+d\Gamma\big((-K_{1}P_{\Delta}K_{1})\oplus
K_{2}P_{\Delta}K_{2}\big).\label{iufoioi}
\end{align}
We finally set
\begin{equation}\label{dsgg}a(\Delta):=a^{+}(\Delta)+a^{0}(\Delta)+a^{-}(\Delta),\quad\Delta\in \mathcal B_0(X),\end{equation} which are linear continuous operators in $\mathcal A\mathcal
F_{\mathrm{fin}}(\mathcal H)$.

\begin{lemma}\label{futyfytd}
The operators $a(\Delta)$, $\Delta\in \mathcal B_0(X)$, commute on
$\mathcal A\mathcal F_{\mathrm{fin}}(\mathcal H)$.
\end{lemma}

\noindent {\it Proof}.  For any $\Delta_1,\Delta_2\in \mathcal B_0(X)$, we trivially have:
\begin{equation}\label{1333}
a^{+}(\Delta_1)a^{+}(\Delta_2)=a^{+}(\Delta_2)a^{+}(\Delta_1).
\end{equation}

Next, we evaluate
\begin{equation}\label{1000}
d\Gamma\big(K_{1}P_{\Delta_1}K_{1}\big)\mathcal N(K_{2}P_{\Delta_2}K_{1})_{2,1}^{\wedge}=\big((\pmb 1\otimes
K_{1}P_{\Delta_1}K_{1})  \mathcal N(K_{2}P_{\Delta_2}K_{1})_{2,1} \big)^{\wedge},
\end{equation}
where $\wedge$ denotes antisymmetrization. For any $u_{1}\in H_{1}$ and $u_{2}\in H_{2}$, we get:
\begin{align*}
&\big((\pmb 1\otimes
K_{1}P_{\Delta_1}K_{1})  \mathcal N(K_{2}P_{\Delta_2}K_{1})_{2,1} ,u_{2}\otimes
u_{1}\big)_{\mathcal H^{\otimes 2}}\\
&\qquad=\big(\mathcal N(K_{2}P_{\Delta_2}K_{1})_{2,1}, u_{2}\otimes K_{1}P_{\Delta_1}K_{1}u_{1}\big)_{\mathcal H^{\otimes 2}}\\
&\qquad=\big(u_{2}, K_{2}P_{\Delta_2}K_{1}K_{1}P_{\Delta_1}K_{1}u_{1}\big)_{H}\\
&\qquad=\big(u_{2}, K_{2}P_{\Delta_2}KP_{\Delta_1}K_{1}u_{1}\big)_{H}.
\end{align*}
Therefore, $(\pmb 1\otimes
K_{1}P_{\Delta_1}K_{1})  \mathcal N(K_{2}P_{\Delta_2}K_{1})_{2,1}$ is the kernel of
the operator $K_{2}P_{\Delta_2}KP_{\Delta_1}K_{1}$ realized as the element of
$H_{2}\otimes H_{1}$. We denote it by $\mathcal
N(K_{2}P_{\Delta_2}KP_{\Delta_1}K_{1})_{2,1}$. Therefore, by (\ref{1000}),
\begin{equation}\label{1111}
d\Gamma\big(K_{1}P_{\Delta_1}K_{1}\big)\mathcal N(K_{2}P_{\Delta_2}K_{1})_{2,1}^{\wedge}=\mathcal
N\big(K_{2}P_{\Delta_2}KP_{\Delta_1}K_{1}\big)_{2,1}^{\wedge}.
\end{equation}

Analogously, we get, for any $u_{1}\in H_{1}$, $u_{2}\in H_{2}$,
$$
\big((K_{2}P_{\Delta_1}K_{2}\otimes \pmb 1) \mathcal N(K_{2}P_{\Delta_2}K_{1})_{2,1}, u_{2}\otimes
u_{1}\big)_{\mathcal H^{\otimes 2}}
=\big(u_{2}, K_{2}P_{\Delta_1}(\pmb 1-K)P_{\Delta_2}K_{1}u_{1}\big)_{H},
$$
and hence,
\begin{equation}\label{1222}
d\Gamma\big(K_{2}P_{\Delta_1}K_{2}\big)\mathcal N(K_{2}P_{\Delta_2}K_{1})_{2,1}^{\wedge}=\mathcal
N(K_{2}P_{\Delta_1}(\pmb 1-K)P_{\Delta_2}K_{1})_{2,1}^{\wedge}.
\end{equation}

By (\ref{1111}) and (\ref{1222}), a straightforward calculation shows
that
\begin{equation}\label{1444}
a^{0}(\Delta_1)a^{+}(\Delta_2)+a^{+}(\Delta_1)a^{0}(\Delta_2)=a^{0}(\Delta_2)a^{+}(\Delta_1)+a^{+}(\Delta_2)a^{0}(\Delta_1).
\end{equation}

Next, by (\ref{777}), we have:
\begin{align}\label{1888}
&a^{-}(\Delta_1)a^{+}(\Delta_2)h_{1}\wedge\dots\wedge h_{n}\nonumber\\
&\qquad=\big((\mathcal N(K_2P_{\Delta_1}K_1),\mathcal N(K_2P_{\Delta_2}K_1))
_{H^{\otimes 2}}\pmb
1
+a^{+}(\Delta_2)a^{-}(\Delta_1)\big)h_{1}\wedge\dots\wedge h_{n}\nonumber\\
&\qquad\quad-\sum_{i=1}^{n}h_{1}\wedge\dots\wedge
h_{i-1}\notag\\
&\qquad \wedge\bigg(\int_{X}\int_{X} \mathcal N(K_2P_{\Delta_2}K_1)
(x,\cdot)h_{i}^{(1)}(y)
 \mathcal N(K_2P_{\Delta_1}K_1)
(x,y)\,\sigma(dx)\,\sigma(dy)\nonumber\\
&\qquad\quad+\int_{X}\int_{X} \mathcal N(K_2P_{\Delta_2}K_1)
(\cdot,y)h_{i}^{(2)}(x)\mathcal
N(K_2P_{\Delta_1}K_1)(x,y)\,\sigma(dx)\,\sigma(dy)\bigg)\notag\\
&\qquad\wedge h_{i+1}\wedge\dots\wedge h_{n}.
\end{align}
For any $u\in H$,
\begin{align}
&\bigg(\int_{X}\int_{X}\mathcal N(K_2P_{\Delta_2}K_1)(x,\cdot)h_{i}^{(1)}(y)
\mathcal N(K_2P_{\Delta_1}K_1)(x,y)\,\sigma(dx)\,\sigma(dy),u\bigg)_{H}\nonumber\\
&\qquad=\int_{X}\int_{X}\int_{X} \mathcal N(K_2P_{\Delta_2}K_1)(x,z)h_{i}^{(1)}
(y)\mathcal N(K_2P_{\Delta_1}K_1)(x,y)u(z)\,\sigma(dx)\,\sigma(dy)\,\sigma(dz)\nonumber\\
&\qquad=\int_{X}\sigma(dy)\int_{X}\sigma(dz)\,h_{i}^{(1)}(y)u(z)\int_{X}\sigma(dx)\,\mathcal
N(K_1P_{\Delta_2}K_2)(z,x)\mathcal N(K_2P_{\Delta_1}K_1)(x,y)\notag\\
&\qquad=\int_{X}\sigma(dy)\int_{X}\sigma(dz)\,h_{i}^{(1)}(y)u(z)\mathcal
N(K_{1}P_{\Delta_2}(\pmb 1-K)P_{\Delta_1}K_{1})(z,y)\notag\\&\qquad =\big(K_{1}P_{\Delta_2}(\pmb
1-K)P_{\Delta_1}K_{1}h_{i}^{(1)},u\big)_H. \notag
\end{align}
 Therefore,
\begin{align}\label{1666}
&\int_{X}\int_{X}\mathcal N(K_2P_{\Delta_2}K_1)(x,\cdot)h_{i}^{(1)}
(y)\mathcal N(K_2P_{\Delta_1}K_1)(x,y)\,\sigma(dx)\,\sigma(dy)\nonumber\\
&\qquad=K_{1}P_{\Delta_2}(\pmb 1-K)P_{\Delta_1}K_{1}h_{i}^{(1)}.
\end{align}
Analogously
\begin{align}\label{1777}
&\int_{X}\int_{X} \mathcal N(K_2P_{\Delta_2}K_1) (\cdot,y)h_{i}^{(2)}(x)\mathcal
N(K_2P_{\Delta_1}K_1)(x,y)\,\sigma(dx)\,\sigma(dy)\nonumber\\
&\qquad=K_{2}P_{\Delta_2}KP_{\Delta_1}K_{2}h_{i}^{(2)}.
\end{align}
By (\ref{1888})--(\ref{1777}),
\begin{align}\label{1999}
&a^{-}(\Delta_1)a^{+}(\Delta_2)h_{1}\wedge\dots\wedge h_{n}\nonumber\\
&\qquad=\big((\mathcal N(K_2P_{\Delta_1}K_1),\mathcal N(K_2P_{\Delta_2}K_1))_{H^{\otimes 2}}\pmb
1
+a^{+}(\Delta_2)a^{-}(\Delta_1)\big)h_{1}\wedge\dots\wedge h_{n}\nonumber\\
&\qquad\quad-d\Gamma\big((K_{1}P_{\Delta_2}(\pmb 1-K)P_{\Delta_1}K_{1})\oplus
(K_{2}P_{\Delta_2}KP_{\Delta_1}K_{2})\big).
\end{align}
Using (\ref{1999}), we conclude that
\begin{align}\label{2000}
&a^{+}(\Delta_1)a^{-}(\Delta_2)+a^{-}(\Delta_1)a^{+}(\Delta_2)+a^{0}(\Delta_1)a^{0}(\Delta_2)\nonumber\\
&\qquad=a^{+}(\Delta_2)a^{-}(\Delta_1)+a^{-}(\Delta_2)a^{+}(\Delta_1)+a^{0}(\Delta_2)a^{0}(\Delta_1).
\end{align}

By (\ref{1333}), (\ref{1444}), (\ref{2000}) and the equalities obtained by
taking the adjoint operators in (\ref{1333}), (\ref{1444}), we conclude the
statement of the lemma.\quad $\square$\vspace{2mm}

We will now show that the family $(a(\Delta))_{\Delta\in\mathcal
B_0(X)}$ has a correlation measure $\rho$ with respect to the vacuum vector $\Omega$.
 Using \eqref{uiyftur}, we informally compute that, for any $\Delta_1\dots,\Delta_n\in\mathcal B_0(X)$,
\begin{equation}\label{uftyukt} \mathcal Q(\chi_{\Delta_1}\odot\dots\odot\chi_{\Delta_n})=\frac{1}{n!}
\int_{\Delta_1}\sigma(dx_1)\dotsm\int_{\Delta_n}\sigma(dx_n)\,
\Psi^*(x_n)\dotsm \Psi^*(x_1)\Psi(x_1)\dotsm\Psi(x_n),\end{equation}
so that
\begin{align*} &Q(\chi_{\Delta_1}\odot\dots\odot\chi_{\Delta_n})\\
&\qquad =\frac{1}{n!}
\int_{\Delta_1}\sigma(dx_1)\dotsm\int_{\Delta_n}\sigma(dx_n)\,
\Psi^*(x_n)\dotsm \Psi^*(x_1)\Phi_{1}^{*}(\varkappa_{1,x_{1}})\dots
\Phi_{1}^{*}(\varkappa_{1,x_{n}})\Omega.\end{align*} Hence,  we need to make
sense out of the following operators, which are heuristically given by
\begin{align}
&\mathcal T(\Delta_1, \dots, \Delta_n) \notag\\
&\quad=\int_{\Delta_n}\sigma(dx_{n})\, \big(\Phi_2^*(\varkappa_{2,x_n})+\Phi_1(\varkappa_{1,x_n})
\big)\notag\\ &\quad\quad\times
 \bigg(\int_{\Delta_{n-1}}\sigma(dx_{n-1})
\,\big(\Phi_2^*(\varkappa_{2,x_{n-1}})+\Phi_1(\varkappa_{1,x_{n-1}})\big)\times\dotsm\notag\\
&\quad\quad\times \bigg(\int_{\Delta_1}\sigma(dx_{1})
\,\big(\Phi_2^*(\varkappa_{2,x_1})+\Phi_1(\varkappa_{1,x_1})
\big)\Phi_{1}^{*}(\varkappa_{1,x_{1}})\bigg)
\dotsm
\Phi_{1}^{*}(\varkappa_{1,x_{n-1}})\bigg)\Phi_{1}^{*}(\varkappa_{1,x_{n}}).
\label{yygyug}
\end{align}

For Hilbert spaces $\mathfrak H_1$ and $\mathfrak H_2$, we denote by $\mathcal L(\mathfrak H_1,\mathfrak H_2)$ the Banach space of linear continuous operators from $\mathfrak H_1$ into $\mathfrak H_2$.
Let $R_{k,n}\in\mathcal L(\mathcal H^{\wedge k}, \mathcal H^{\wedge n})$ and let $\Delta\in\mathcal B_0(X)$.

Taking into account \eqref{444} and \eqref{555}, we define an operator
$$
\int_\Delta\sigma(dx)\, \Phi^*_2(\varkappa_{2,x})R_{k,n}\Phi^*(\varkappa_{1,x})
\in\mathcal L(\mathcal{H}^{\wedge(k-1)}, \mathcal{H}^{\wedge(n+1)}) $$
as follows: for each $f^{(k-1)}\in \mathcal{H}^{\wedge(k-1)}$ we set
$$ \int_\Delta\sigma(dx)\, \Phi^*_2(\varkappa_{2,x})R_{k,n}\Phi^*(\varkappa_{1,x})
f^{(k-1)}:=\mathcal P_{n+1}(\pmb 1\otimes (R_{k,n}\mathcal P_k))
(\mathcal N(K_2P_\Delta K_1)\otimes f^{(k-1)}).
$$ Here, $\mathcal P_i$ denotes the orthogonal projection of $\mathcal H^{\otimes i}$ onto $\mathcal H^{\wedge i}$.

Next, using \eqref{*}, we easily conclude that the operator-valued function
$$ X\ni x\mapsto \Phi_1(\varkappa_{1,x})R_{k,n}
\Phi_1^*(\varkappa_{1,x})\in\mathcal L(\mathcal H^{\wedge(k-1)},\mathcal H^{\wedge(n-1)})$$
is strongly measurable, and Bochner-integrable over $\Delta$ (see e.g.\ \cite{BUS} for details on Bochner integral).
So, we define $$ \int_\Delta \sigma(dx)\, \Phi_1(\varkappa_{1,x})R_{k,n}
\Phi_1^*(\varkappa_{1,x})\in\mathcal L(\mathcal H^{\wedge(k-1)},\mathcal H^{\wedge(n-1)}) $$ as a Bochner integral.

Finally, by linearity, for any linear continuous operator $R$ in $\mathcal{AF_{\text{fin}}}(\mathcal H)$, we define
$\int_\Delta\sigma(dx)\, \Phi^*_2(\varkappa_{2,x})R\Phi^*(\varkappa_{1,x})$
and $\int_\Delta \sigma(dx)\, \Phi_1(\varkappa_{1,x})R
\Phi_1^*(\varkappa_{1,x})$
as linear continuous operators  in $\mathcal{AF_{\text{fin}}}(\mathcal H)$.
Hence, by induction, the operator \eqref{yygyug} is well defined.

\begin{lemma}\label{dtstes}
For each $n\in\mathbb N$ and any $\Delta_1,\dots,\Delta_n\in\mathcal B_0(X)$, we
have\rom: $$
Q(\chi_{\Delta_1}\odot\dots\odot\chi_{\Delta_n})=\frac1{n!}\,\mathcal T
(\Delta_1, \dots, \Delta_n)\Omega.$$
\end{lemma}

\noindent {\it Proof.}  We first state that, for any $\Delta_1,\Delta_2\in\mathcal B_0(X)$ and any linear continuous operator
$R$ in $\mathcal{AF}_{\mathrm{fin}}(\mathcal H)$, we have
\begin{align} &a(\Delta_1)\int_{\Delta_2}\sigma(dx)\,(\Phi^*_2(\varkappa_{2,x})
+\Phi_1(\varkappa_{1,x}))R\Phi^*(\varkappa_{1,x})\notag\\
&\qquad=\int_{\Delta_2}\sigma(dx)\,(\Phi^*_2(\varkappa_{2,x})
+\Phi_1(\varkappa_{1,x}))a(\Delta_1)R\Phi^*(\varkappa_{1,x})\notag\\
&\qquad\quad-\int_{\Delta_1\cap\Delta_2}\sigma(dx)\,(\Phi^*_2(\varkappa_{2,x})
+\Phi_1(\varkappa_{1,x}))R\Phi^*(\varkappa_{1,x}).\label{jiuuu}
\end{align}
 Indeed, to show \eqref{jiuuu} it is sufficient to consider
 the case where $R=R_{k,n}\in\mathcal L(\mathcal H^{\wedge k},\mathcal H^{\wedge n})$
  and $R_{n,k}$ has the form
 $$ R_{k,n}f^{(k)}=(f^{(k)},u_1\wedge\dots\wedge u_k)_{\mathcal H^{\wedge k}}\,v_1\wedge\dots\wedge v_n,\qquad
 f^{(k)}\in\mathcal H^{\wedge k},$$
 with $u_1,\dots,u_k,v_1,\dots,v_n\in \mathcal H$. But   \eqref{jiuuu}  with $R=R_{k,n}$
  of such a form can be deduced analogously to the proof of Lemma~\ref{futyfytd}. Now, by virtue of  the recurrence formula \eqref{uiifdyt}, the statement
 of Lemma~\ref{dtstes}  follows from \eqref{jiuuu} by induction. \quad $\square$

\begin{remark}\rom{It is, in fact, possible to rigorously define the operator on the right hand side of \eqref{uftyukt}, and show that equality \eqref{uftyukt} indeed holds.
}\end{remark}

\begin{lemma}\label{iuqd}
The family of operators $(a(\Delta))_{\Delta\in\mathcal B_0(X)}$ has a correlation measure $\rho$ with respect to $\Omega$, and the  restriction
of $\rho$ to $(\Gamma_X^{(n)},\mathcal B(\Gamma_X^{(n)}))$ is given by
\begin{equation}\label{hyfuf} \rho^{(n)}(dx_1,\dots,dx_n)=\frac1{n!}\,\operatorname{det}[k(x_i,x_j)]_{i,j=1}^n\,\sigma(dx_1)\dotsm
\sigma(dx_n)\end{equation}
\rom(recall that we have identified $\mathcal B(\ddot\Gamma_X^{(n)})$ with $\mathcal{B}_{\mathrm{sym}}(X^{n})$, and 
$\mathcal B(\Gamma_X^{(n)})\subset \mathcal B(\ddot\Gamma_X^{(n)})$\rom).
\end{lemma}

\noindent {\it Proof}. By  \eqref{yygyug} and Lemma \ref{dtstes}, for each $n\in\mathbb N$ and any $\Delta_1,\dots,\Delta_n
\in\mathcal B_0(X)$, we have
\begin{align}
&(Q(\chi_{\Delta_1}\odot\dots\odot\chi_{\Delta_n}),\Omega)_{\mathcal{AF}(\mathcal H)}\notag\\
&\qquad=\frac1{n!}\bigg( \int_{\Delta_n}\sigma(dx_{n})\, \Phi_1(\varkappa_{1,x_n})
\bigg(\int_{\Delta_{n-1}}\sigma(dx_{n-1})
\, \Phi_1(\varkappa_{1,x_{n-1}})\notag\\
&\qquad \quad \times\dotsm \times \bigg(\int_{\Delta_1}\sigma(dx_{1})
\, \Phi_1(\varkappa_{1,x_1})
\Phi_{1}^{*}(\varkappa_{1,x_{1}})\bigg)\dotsm
\Phi_{1}^{*}(\varkappa_{1,x_{n-1}})\bigg)\Phi_{1}^{*}(\varkappa_{1,x_{n}})\Omega,\Omega\bigg)_{\mathcal{AF}(\mathcal H)}
\notag\\
&\qquad =\int_{\Delta_n}\sigma(dx_n)\dotsm\int_{\Delta_1}\sigma(dx_1)\big\| \varkappa_{1,x_1}\wedge\dots\wedge\varkappa_{1,x_n}
\big\|^2 _{\mathcal H^{\wedge n}}\label{frdtrsd}\\
&\qquad =\frac1{n!}\int_{\Delta_n}\sigma(dx_n)\dotsm\int_{\Delta_1}\sigma(dx_1)\operatorname{det}[k(x_i,x_j)]_{i,j=1}^n\,\sigma(dx_1)\dotsm
\sigma(dx_n).\label{uyfr}
\end{align}
Note that, by \eqref{frdtrsd}, the right hand side of \eqref{hyfuf}  indeed defines a measure. Hence, the statement
of the lemma follows from \eqref{uyfr}\quad $\square$

\begin{lemma}\label{uig}
The correlation measure given in \eqref{hyfuf} satisfies \rom{(LB)}.
\end{lemma}

\noindent{\it Proof}. For each $\Delta\in\mathcal B_0(X)$ and $n\in\mathbb N$, we evidently have
\begin{align*}\rho(\Gamma_\Delta^{(n)})&\le \bigg(\int_\Delta \|\varkappa_{1,x}\|_H^2\,\sigma(dx)\bigg)^n\\
&=\bigg(\int_{\Delta}\int_{X}\mathcal N(K_1)(x,y)^2\,\sigma(dx)\,\sigma(dy)\bigg)^n,
\end{align*}
from where the statement follows.\quad $\square$

By Lemmas \ref{futyfytd}, \ref{iuqd}, \ref{uig}, and  Theorem \ref{uytwa}, we get

\begin{theorem}\label{hyufgu}
For the family $(a(\Delta))_{\Delta\in\mathcal B_0(X)}$ defined
 through formulas \eqref{jhvsu},
\eqref{777}, \eqref{iufoioi}, and \eqref{dsgg}, the statement of
Theorem~\rom{\ref{uytwa}} holds with the correlation measure given by
\eqref{hyfuf}.
\end{theorem}

Let us now briefly mention the boson case. About the operator $K$ we make the
same assumptions as in the fermion case, apart from the assumption that
$K\le\pmb1 $. We  set $K_1:=\sqrt K$ (just as above) and
$K_2:=(\pmb1+K)^{1/2}$. We then essentially repeat the fermion case, using
however the symmetric Fock space $\mathcal {SF}(\mathcal H)$ instead of the
antysymmetric Fock space $\mathcal {AF}(\mathcal H)$. The operators $\Psi(f)$,
$\Psi^*(f)$ (see \eqref{vfsiu}) now satisfy the CCR (use the commutator
$[A,B]_-:=AB-BA$ instead of the anticommutator in \eqref{uiyftur}). The
counterpart of formulas \eqref{frdtrsd}, \eqref{uyfr} reads as folllows:
\begin{align*}
&(Q(\chi_{\Delta_1}\odot\dots\odot\chi_{\Delta_n}),\Omega)_{\mathcal{AF}(\mathcal H)}\notag\\
&\qquad =\int_{\Delta_n}\sigma(dx_n)\dotsm\int_{\Delta_1}\sigma(dx_1)\big\| \varkappa_{1,x_1}\odot\dots\odot\varkappa_{1,x_n}
\big\|^2 _{\mathcal H^{\odot n}}\\
&\qquad =\frac1{n!}\int_{\Delta_n}\sigma(dx_n)\dotsm\int_{\Delta_1}\sigma(dx_1)\operatorname{per}[k(x_i,x_j)]_{i,j=1,\dots,n}\,\sigma(dx_1)\dotsm
\sigma(dx_n).
\end{align*}
Thus the corresponding correlation measure is given by \eqref{hyfuf} in which the determinant is replaced by the permanent.

\section{Reduced particle densities}\label{yuyuu}
Let the operators $K$, $K_1$, $K_2$ be as in the fermion part of Section~\ref{iusfdgh}. Let $l\in\mathbb N$,
$l\ge2$, and we now take $2l$ copies of the Hilbert space $H=L^2(X,\sigma)$: $H_{1,i}$ and $H_{2,i}$, $i=1,\dots,l$. We denote
$\mathcal H^{(l)}:=\bigoplus_{i=1}^l(H_{1,i}\oplus H_{2,i})$.
For each $f\in H$, we consider the following operators in $\mathcal{AF}(\mathcal H^{(l)})$:
\begin{align}
\Psi(f):=&\sum_{i=1}^l\big(\Phi_{2,i}(l^{-1/2}K_{2}f)+\Phi_{1,i}^{*}(l^{-1/2}K_{1}f)\big),\notag\\
\Psi^{*}(f):=&\sum_{i=1}^l \big(\Phi_{2,i}^{*}(l^{-1/2}K_{2}f)+\Phi_{1,i}(l^{-1/2}K_{1}f)\big)\label{juigyui}
\end{align}
(we are using obvious notations, analogous to those of Section~\ref{iusfdgh}).
It is easy to see that these operators  satisfy the CAR \eqref{uiyftur}. Furthermore,
the $n$-point functions of this representation of the CAR are again given by \eqref{uytfrtyrfrr}.
Therefore, the representation of the CAR given by \eqref{juigyui} is unitary equivalent to the representation
  \eqref{vfsiu}.

  The particle density of the representation \eqref{juigyui} is heuristically given by
  \begin{align}
  a^{(l)}(x):& =\Psi^*(x) \Psi(x)\notag\\
  &=\sum_{i=1}^l\sum_{j=1}^l
\big(\Phi_{2,i}^*(l^{-1/2}\varkappa_{2,i,x})+\Phi_{1,i}(l^{-1/2}\varkappa_{1,i,x})\big)\notag\\
&\qquad\quad\times  \big(\Phi_{2,j}(l^{-1/2}\varkappa_{2,i,x})+\Phi_{1,j}^*(l^{-1/2}\varkappa_{1,i,x})\big).\label{yur}
\end{align}
One can rigorously construct a corresponding family of commuting Hermitian operators,
$(a^{(l)}(\Delta))_{\Delta\in\mathcal B_0(X)}$, and show that, as expected, the family $(a^{(l)}(\Delta))_{\Delta\in\mathcal B_0(X)}$ has the same correlation measure
\eqref{hyfuf} with respect to the vacuum vector $\Omega$.

Now, let us consider the reduced  particle density
\begin{equation} \label{jukiguyg}R^{(l)}(x) :=\sum_{i=1}^l
\big(\Phi_{2,i}^*(l^{-1/2}\varkappa_{2,i,x})+\Phi_{1,i}(l^{-1/2}\varkappa_{1,i,x})\big)  \big(\Phi_{2,i}(l^{-1/2}\varkappa_{2,i,x})+\Phi_{1,j}^*(l^{-1/2}\varkappa_{1,i,x})\big)\end{equation}
(i.e., we have taken only the ``diagonal elements'' of the double sum).
Analogously to Section~\ref{iusfdgh}, one can rigorously realize $(R^{(l)}(\Delta))_{\Delta\in\mathcal B_0(X)}$
as a  family of commuting Hermitian operators in $\mathcal {AF}(\mathcal H^{(l)})$.
The counterpart of formulas
\eqref{frdtrsd}, \eqref{uyfr} reads as folllows:
\begin{align}
&(Q(\chi_{\Delta_1}\odot\dots\odot\chi_{\Delta_n}),\Omega)_{\mathcal{AF}(\mathcal H^{(l)})}\notag\\
&=\sum_{i_1=1}^l\dots\sum_{i_n=1}^l \frac1{n!}\bigg( \int_{\Delta_n}\sigma(dx_{n})\, \Phi_{1,i_n}(l^{-1/2}\varkappa_{1,i_n,x_n})
\notag\\
& \quad \times\dotsm \times \bigg(\int_{\Delta_1}\sigma(dx_{1})
\,\Phi_{1,i_1}(l^{-1/2}\varkappa_{1,i_1,x_1})
\Phi_{1,i_1}^{*}(l^{-1/2}\varkappa_{1,i_1,x_{1}}) \bigg)\notag\\
&\quad\dotsm
\Phi_{1,i_n}^{*}(l^{-1/2}\varkappa_{1,i_n,x_{n}}) \Omega,\Omega\bigg)_{\mathcal{AF}(\mathcal H^{(l)})}
\notag\\
& =\int_{\Delta_n}\sigma(dx_n)\dotsm\int_{\Delta_1}\sigma(dx_1)\, \sum_{i_1=1}^l\dots\sum_{i_n=1}^l
\big\| (l^{-1/2}\varkappa_{1,i_1,x_1}) \wedge\dots\wedge (l^{-1/2}\varkappa_{1,i_n,x_n} ))\big \|^2_{(\mathcal H^{(l)})^{\wedge n}}.\label{jfcdydc}
\end{align}
Hence, $(R^{(l)}(\Delta))_{\Delta\in\mathcal B_0(X)}$ has the correlation measure, whose restriction to $(\Gamma_X^{(n)},\mathcal B(\Gamma_X^{(n)}))$ is given by
\begin{multline*}\rho^{(n)}(dx_1,\dots,dx_n)\\=\sum_{i_1=1}^l\dots\sum_{i_n=1}^l
\| (l^{-1/2}\varkappa_{1,i_1,x_1}) \wedge\dots\wedge (l^{-1/2}\varkappa_{1,i_n,x_n} )) \|^2_{(\mathcal H^{(l)})^{\wedge n}}\, \sigma(dx_1)\dotsm
\sigma(dx_n).\end{multline*}
A combinatoric exercise shows that the $\rho^{(n)}$ can be written in the form
\begin{equation}\label{ig} \rho^{(n)}(dx_1,\dots,dx_n)=\frac1{n!}\,\operatorname{det}_{-1/l}[l\times k(x_i,x_j)]_{i,j=1}^n\,\sigma(dx_1)\dotsm
\sigma(dx_n).\end{equation}
Here, for any $\alpha\in\mathbb R$ and a square matrix $A=(a_{i,j})_{i,j=1}^n$, $\operatorname{det}_\alpha A$
denotes the Vere-Jones $\alpha$-determinant (see \cite{ST}):
$$ \operatorname{det}_\alpha A:= \sum_{\xi\in S_n}\alpha^{n-\nu(\xi)}\prod_{i=1}^n a_{i,\xi(i)},$$
where $\nu(\xi)$ denotes the number of cycles in the permutation $\xi$.

The correlation measure \eqref{ig} satisfies (LB), and so the statement of Theorem~\ref{uytwa} holds
for the family $(R^{(l)}(\Delta))_{\Delta\in\mathcal B_0(X)}$. Formula \eqref{ig} also shows that the corresponding measure on $\Gamma_X$, which we denote by $\mu^{(l)}$ is the fermion-like point process considered in \cite{ST}.

It is heuristically clear from \eqref{jukiguyg} that the measure $\mu^{(l)}$ is the $l$-fold convolution of
fermion point processes  corresponding to the operator $K/l$. This, in fact, can be  rigorously shown, since the correlation
measure of a convolution of point processes may be easily expressed in terms of the correlation measures
of the initial point processes, see also \cite{ST}.

Finally, an analogous construction can be  carried out in the boson case, leading
to the correlation function \eqref{ig}, in which $\operatorname{det}_{-1/l}$ is replaced by $\operatorname{det}_{1/l}$.

 \begin{center}
{\bf Acknowledgements}\end{center}

 The first named author acknowledges the financial support of the SFB 701 `` Spectral
structures and topological methods in mathematics'', Bielefeld University.
We are grateful to Yuri Kondratiev and Tobias Kuna  for many useful discussions.

\end{document}